\documentclass{amsart}

\usepackage{amsfonts}
\usepackage{amsmath}
\usepackage{amsthm}
\usepackage{amssymb}
\usepackage[english]{babel}
\usepackage{graphics}
\usepackage{latexsym}
\usepackage{longtable} \usepackage{mathrsfs}
\usepackage{graphicx}
\usepackage{wrapfig} 
\usepackage[pdftex,colorlinks=true]{hyperref}

\newtheorem{theorem}{Theorem}
\newtheorem{corollary}[theorem]{Corollary}
\newtheorem{lemma}[theorem]{Lemma}

\newtheorem{claim}[theorem]{Claim}
\newtheorem{example}[theorem]{Example}

\theoremstyle{definition}
\newtheorem{definition}[theorem]{Definition}

\renewcommand{\S}{\mathcal{S}}

\newcommand{\D}{\mathrm{D}}
\newcommand{\A}{\textbf{A}}

\newcommand{\R}{\mathbb{R}}

\newcommand{\noi}{\noindent}
\newcommand{\ms}{\medskip}

\newcommand{\al}{\alpha}
\newcommand{\be}{\beta}
\newcommand{\ga}{\gamma}

\newcommand{\e}{\varepsilon}

\newcommand{\la}{\lambda}
\newcommand{\ka}{\kappa}
\newcommand{\Om}{\Omega}

\newcommand{\larrow}{\longrightarrow}
\newcommand{\ot}{\otimes}

\newcommand{\ri}{\rightarrow}

\newcommand{\p}{\partial}
\newcommand{\sub}{\subseteq}

\newcommand{\by}{\times}

\newcommand{\sgn}{\mathrm{sgn}}
\newcommand{\ess}{\mathrm{ess}}

\newcommand{\cof}{\mathrm{cof}}

\newcommand{\bt}{\begin{theorem}}\newcommand{\et}{\end{theorem}}
\newcommand{\bd}{\begin{definition}}\newcommand{\ed}{\end{definition}}
\newcommand{\bl}{\begin{lemma}}\newcommand{\el}{\end{lemma}}
\newcommand{\beq}{\begin{equation}}\newcommand{\eeq}{\end{equation}}
\newcommand{\bc}{\begin{claim}}\newcommand{\ec}{\end{claim}}
\newcommand{\bex}{\begin{example}}\newcommand{\eex}{\end{example}}
\newcommand{\bcor}{\begin{corollary}}\newcommand{\ecor}{\end{corollary}}
\newcommand{\bp}{\begin{proof}}\newcommand{\ep}{\end{proof}}

\newcommand{\BPL}{\medskip \noindent \textbf{Proof of Lemma} }

\newcommand{\BPT}{\medskip \noindent \textbf{Proof of Theorem} }

\numberwithin{equation}{section}

\begin{document}

\title[Well-Posedness of Global Fully Nonlinear Systems]{On the Well-Posedness of Global Fully Nonlinear First Order Elliptic Systems}

\author{Hussien Abugirda}
\address{Department of Mathematics, College of Science, University of Basra, Basra, Iraq \& Department of Mathematics and Statistics, University of Reading, Whiteknights, PO Box 220, Reading RG6 6AX, UK}
\email{h.a.h.abugirda@student.reading.ac.uk}

\author{Nikos Katzourakis}
\address{Department of Mathematics and Statistics, University of Reading, Whiteknights, PO Box 220, Reading RG6 6AX, UK}
\email{ n.katzourakis@reading.ac.uk}


\date{}



\begin{abstract} In the very recent paper \cite{K1}, the second author proved that for any $ f\in L^2(\R^n,\R^N)$, the fully nonlinear first order system $F(\cdot,\mathrm{D} u) =f$ is well posed in the so-called J.L.\ Lions space and moreover the unique strong solution $u:\R^n\larrow \R^N$ to the problem satisfies a quantitative estimate. A central ingredient in the proof was the introduction of an appropriate notion of ellipticity for $F$ inspired by Campanato's classical work in the 2nd order case. Herein we extend the results of \cite{K1} by introducing a new strictly weaker ellipticity condition and by proving well posedness in the same ``energy" space.
\end{abstract}

\maketitle

\section{Introduction} \label{section1}

In this paper we consider the problem of existence and uniqueness of global strong solutions $u:\R^n \larrow \R^N$ to the fully nonlinear first order PDE system
\beq  \label{1.1}
F(\cdot,\D u) \,=\, f, \ \ \text{ a.e.\ on }\R^n,
\eeq
where $n,N\geq 2$ and $F : \R^{n} \by \R^{Nn} \larrow  \R^N$ is a  Carath\'eodory map. The latter means that $F(\cdot,X)$ is a measurable map for all $X \in \R^{Nn}$ and  $F(x,\cdot)$ is a  continuous map for almost every $x \in \R^{n}$. The gradient $\D u :\R^n \larrow \R^{Nn}$ of our solution $u=(u_1,...,u_N)^\top$ is viewed as an $N\by n$ matrix-valued map $\D u=(\D_i u_\al)_{i=1...n}^{\al=1...N}$ and the right hand side $f$ is assumed to be in $L^2(\R^n,\R^N)$.

The method we use in this paper to study \eqref{1.1} follows that of the recent paper \cite{K1} of the second author. Therein the author introduced and employed a new perturbation method in order to solve \eqref{1.1} which is based on the solvability of the respective linearised system and a structural ellipticity hypothesis on $F$, inspired by the classical work of Campanato in the fully nonlinear second order case $\mathcal{F}(\cdot,\D^2u)=f$ (see \cite{C0}-\cite{C5}, \cite{Co1,Co2}  and \cite{Ta1}-\cite{Ta3}). Loosely speaking, the ellipticity notion of \cite{K1} requires that $F$ is ``not too far away" from a linear constant coefficient first order differential operator. In the linear case of constant coefficients, $F$ assumes the form
\[
F(x,X)\, =\, \sum_{\al,\be=1}^N \sum_{j=1}^n\A_{\al \be j} X_{\be j}\,e^\al,
\]
for some linear map $\A : \R^{Nn}\larrow \R^N$. We will follow almost the same conventions as in \cite{K1}, for instance we will denote the standard bases of $\R^n$, $\R^N$ and $\R^{N\by n}$ by  $\{e^i\}$, $\{e^\al\}$ and $\{e^\al \ot e^i\}$ respectively. In the linear case, \eqref{1.1} can be written as
\[
\sum_{\be=1}^N \sum_{j=1}^n\A_{\al \be j} \D_ju_\be \,=\, f_\al,\ \ \ \ \al=1,...,N,
\]
and compactly in vector notation as
\beq \label{1.2}
\A: \D u \,=\, f.
\eeq
The appropriate well-known notion of ellipticity in the linear case is that the \emph{nullspace of the linear map $\A$ contains no rank-one lines}. This requirement can be quantified as 
\beq\label{1.3}
|\A :\xi \ot a | \, >\, 0,\ \text{ when } \ \xi \neq 0,  \ a \neq 0
\eeq 
which says that all rank-one directions $\xi \ot a \in \R^{Nn}$ are transversal to the nullspace. A prototypical example of such operator $\A : \R^{2 \by 2} \larrow \R^2$ is given by
\beq \label{CR}
\A \, =\, 
\left[
\begin{array}{cc|cc}
1 & 0 & 0 & 1\\
0 & \!\!\!-1 & 1 & 0 
\end{array}
\right]
\eeq
and corresponds to the Cauchy-Riemann PDEs. In \cite{K1} the system \eqref{1.1} was proved to be well-posed by solving \eqref{1.2} via Fourier transform methods and by utilising the following ellipticity notion: \eqref{1.1} is an elliptic system (or $F$ is elliptic) when there exists a linear map
\[
\A \ : \ \R^{Nn} \larrow \R^N
\]
which is elliptic in the sense of \eqref{1.3} and 
\beq \label{1.4}
\underset{x\in \R^n}{\ess\,\sup} \sup_{X,Y\in \R^{Nn},X\neq Y}\frac{\left| \big[F(x,Y) -F(x,X)\big] -\A:(Y-X)\right| }{|Y-X|}\, <\, \nu (A),
\eeq
 where 
\beq  \label{1.5}
 \nu(\A)\ :=\ \min_{|\eta|=|a|=1} \big|\A : \eta \ot a \big|
\eeq
is the ``ellipticity constant" of $\A$. This notion was called ``K-Condition" in \cite{K1}. The functional space in which well posedness was obtained is the so-called J.L.\ Lions space
\beq \label{1.6}
W^{1;2^*\!,2}(\R^n,\R^N)\, :=\, \Big\{ u\in L^{2^*}(\R^n,\R^N)\ : \ \D u\in L^{2}(\R^n,\R^{Nn})\Big\} .
\eeq
Here $2^*$ is the conjugate Sobolev exponent 
\[
2^*\, =\, \frac{2n}{n-2}
\]
(note that ``$L^{2^*}$" means ``$L^{p}$ for $p={2^*}$", not duality) and the natural norm of the space is
\[
 \| u\|_{W^{1;2^*,2}(\R^n)}\, :=\, \|u\|_{L^{2^*}(\R^n)}\, + \, \|\D u\|_{L^{2}(\R^n)}.
\]
In \cite{K1} only global strong a.e.\ solutions on the whole space were considered and for dimensions $n\geq 3$ and $N\geq 2$, in order to avoid the compatibility difficulties which arise in the case of the Dirichlet problem for first order systems on bounded domains and because the case $n=2$ has been studied quite extensively.

In this paper we follow the method introduced in \cite{K1} and we prove well-posedness of \eqref{1.1} in the space \eqref{1.6} for the same dimensions $n\geq 3$ and $N\geq 2$. This is the content of our Theorem \ref{th2}, whilst we also obtain an a priori quantitative estimate in the form of a ``comparison principle" for the distance of two solutions in terms of the distance of the respective right hand sides of \eqref{1.1}. \textit{The main advance in this paper which distinguishes it from the results obtained in \cite{K1} is that we introduce a new notion of ellipticity for \eqref{1.1} which is strictly weaker than \eqref{1.4}, allowing for more general nonlinearities $F$ to be considered.} Our new hypothesis of ellipticity is inspired by an other recent work of the second author \cite{K2} on the second order case. We will refer to our condition as the ``AK-Condition"  (Definition \ref{def2}). In Examples \ref{ex9}, \ref{ex8} we demonstrate that the new condition is genuinely weaker and hence our results indeed generalise those of \cite{K1}. Further, motivated by \cite{K2} we also introduce a related notion which we call pseudo-monotonicity and examine their connection (Lemma \ref{pr1}). The idea of the proof of our main result Theorem \ref{th2} is based, as in \cite{K1}, on the solvability of the linear system, our ellipticity assumption and on a fixed point argument in the form of  Campanato's near operators, which we recall later for the convenience of the reader (Theorem \ref{th3}). 

We conclude this introduction with some comments which contextualise the standing of the topic and connect to previous contributions by other authors. Linear elliptic PDE systems of the first order are of paramount importance in several branches of Analysis like for instance in Complex and Harmonic Analysis. Therefore, they have been extensively studied in several contexts (see e.g.\ Buchanan-Gilbert \cite{BG}, Begehr-Wen \cite{BW}), including regularity theory of PDE (see chapter 7 of Morrey's exposition \cite{Mo} of the Agmon-Douglis-Nirenberg theory), Differential Inclusions and Compensated Compactness theory (Di Perna \cite{DP}, M\"uller \cite{Mu}),  as well as Geometric Analysis and the theory of differential forms (Csat\'o-Dacorogna-Kneuss \cite{CDK}).  

However, except for the paper \cite{K1} the fully nonlinear system \eqref{1.1} is much less studied and understood. By using the Baire category method of the Dacorogna-Marcellini \cite{DM} (which is the analytic counterpart of Gromov's geometric method of Convex Integration), it can be shown that the Dirichlet problem
\beq \label{1.7}
\left\{
\begin{array}{rl}
F(\cdot,\D u)\,=\, f,  &\ \text{ in }\Om, \smallskip\\
u\,=\, g, & \  \text{ on }\p \Om,
\end{array}
\right.
\eeq
has \emph{infinitely many} strong a.e.\ solutions in $W^{1,\infty}(\Om,\R^N)$, for $\Om \sub \R^n$,  $g$ a Lipschitz map and under certain structural \emph{coercivity} and compatibility assumptions. However, roughly speaking ellipticity and coercivity of $F$ are mutually exclusive. In particular, it is well know that the Dirichlet problem \eqref{1.7} is not well posed when $F$ is either linear or elliptic. 

Further, it is well known that single equations, let alone systems of PDE, in general do not have \emph{classical} solutions. In the scalar case $N=1$, the theory of Viscosity Solutions of Crandall-Ishii-Lions (we refer to \cite{K0} for a pedagogical introduction of the topic) furnishes a very successful setting of \emph{generalised} solutions in which Hamilton-Jacobi PDE enjoy strong existence-uniqueness theorems. However, there is no counterpart of this essentially scalar theory for (non-diagonal) systems. The general approach of this paper is inspired by the classical work of Campanato quoted earlier and in a nutshell consists of imposing an appropriate condition that allows to prove well-posedness in the setting of the intermediate theory of \emph{strong a.e.}\ solutions. Notwithstanding, very recently the second author in \cite{K3} has proposed a new theory of generalised solutions in the context of which he has already obtained existence and uniqueness theorems for second order degenerate elliptic systems. We leave the study of the present problem in the context of ``$\mathcal{D}$-solutions" introduced in \cite{K3} for future work.

\section{Preliminaries} \label{section2}

In this section we collect some results taken from our references which are needed for the main results of this paper. The first one below concerns the existence and uniqueness of solutions to the linear first order system with constant coefficient 
\[
\A : \D u \,=\, f, \ \ \text{ a.e.\ on }\R^n,
\]
 with $\A : \R^{Nn} \larrow \R^N$  elliptic in the sense of \eqref{1.3}, namely when the nullspace of $\A$ does not contain rank-one lines. By the compactness of the torus, it can be rewritten equivalently as
\beq \label{2.1}
|\A : \xi \ot a | \,\geq \, \nu\, |\xi| |a|, \ \ \xi \in \R^N, \, a\in \R^n,
\eeq
for some constant $\nu>0$, which can be chosen to be  the \emph{ellipticity constant of $\A$} given by \eqref{1.5}. One can easily see that \eqref{2.1} can be rephrased as
\beq \label{2.3}
\min_{|a|=1}\big| \det(\A a) \big| \, >\, 0,
\eeq
where $\A a$ is the $N \by N$ matrix given by
\[
\A a \, :=\, \sum_{\al,\be=1}^N \sum_{j=1}^n(\A_{\al \be j}\, a_j )\, e^\al \ot e^\be.
\]
It is easy to exhibit examples of tensors $\A$ satisfying \eqref{2.1}. A map $\A : \R^{2 \by 2} \larrow \R^2$ satisfying it is 
\[
\A \, =\, 
\left[
\begin{array}{cc|cc}
\ka & 0 & 0 & \la\\
0 & \!\!\!-\mu & \nu & 0 
\end{array}
\right],
\]
where $\ka,\la,\mu, \nu>0$. A higher dimensional example of map $\A : \R^{4\by 3}\larrow \R^4$ is 
\[
\A \, =\, 
\left[
\begin{array}{rrr|rrr|rrr|rrr}
1 & 0 & 0    &      0 & \!\!\!-1 & 0    &    0 & 0 & \!\!\!-1   & 0& 0 & 0  \\
0 & 1 & 0    &      1 & 0 & 0    &     0 & 0 & 0   &  0& 0 & \!\!\!-1  \\
0 & 0 & 1    &      0 & 0 & 0    &    1 & 0 & 0   &   0& 1 & 0  \\
0 & 0 & 0    &      0 & 0 & 1    &    0 & \!\!\!-1 &0   &   1& 0 & 0  \\
\end{array}
\right]
\]
which corresponds to the electron equation of Dirac in the case where is no external force. For more details we refer to \cite{K1}. 

\bt[Existence-Uniqueness-Representation, cf.\ \cite{K1}] \label{th1} Let $n\geq 3$, $N\geq 2$,  $\A : \R^{Nn}\larrow \R^N$ a linear map satisfying \eqref{2.1} and $f\in L^2(\R^n,\R^N)$. Then, the system
\[
\A : \D u \,=\, f, \ \ \text{ a.e.\ on }\R^n,
\]
has a unique solution $u$ in the space $W^{1;2^*\!,2}(\R^n,\R^N)$ (see \eqref{1.6}), which also satisfies the estimate
\beq \label{2.4}
 \| u\|_{W^{1;2^*,2}(\R^n)}\, \leq\, C \|f\|_{L^{2}(\R^n)}
\eeq
for some $C>0$ depending only on $\emph{\A}$. Moreover, the solution can be represented explicitly as:
\beq \label{2.5}
u\, =\, -\frac{1}{2\pi i} \lim_{m\ri \infty}
\left\{  \widehat{h_m} \ast
\left[ 
\frac{\ \cof \,(\A \sgn )^\top}{ \det (\A \sgn) \ } \overset{\vee}{f}
\right]^{\wedge} 
\right\}.
\eeq
In \eqref{2.5}, $(h_m)^\infty_1$ is any sequence of even functions in the Schwartz class $\S(\R^n)$ satisfying
\[
\text{$0\, \leq\, h_m(x) \,\leq\, \frac{1}{|x|}$ \ \, and \, \ $h_m(x) \larrow \frac{1}{|x|}$, \ for a.e. $x\in \R^n$,\ \,  as $m\ri \infty$.} 
\]
The limit in \eqref{2.5} is meant in the weak $L^{2^*}$ sense as well as a.e.\ on $\R^n$, and $u$ is independent of the choice of sequence $(h_m)^\infty_1 $.
\et

In the above statement,  ``sgn", ``cof" and ``det" symbolise the sign function on $\R^n$, the cofactor and the determinant on $\R^{N \by N}$ respectively. Although the formula \eqref{2.5} involves complex quantities, $u$ above is a \emph{real} vectorial solution. Moreover, the symbol ``$ \widehat{\phantom{...}}$" stands for Fourier transform (with the conventions of \cite{F}) and ``$\overset{\vee}{\phantom{...}}$" stands for its inverse.

\ms

Next, we recall the strict ellipticity condition of the second author taken from \cite{K1} in an alternative form which is more convenient for our analysis. We will relax it in the next section. Let 
\[
\A \ : \ \R^{Nn} \larrow \R^N
\]
be a fixed reference linear map satisfying \eqref{2.1}. 

\begin{definition}[K-Condition of ellipticity, cf.\ \cite{K1}] \label{def1}

Let $F : \R^{n} \by \R^{N n} \larrow  \R^N$ be a Carath\'eodory map. We say  that $F$ is elliptic with respect to $\A$ when there exists $0<\be<1$ such that for all $ X,Y \in \R^{Nn} $ and a.e.\ $x\in \R^n$, we have
\beq \label{3.3}
\Big| \Big[F(x,X+Y) -F(x,X)\Big] -\A:Y\Big| \, \leq\, \be\, \nu (\A) \, {|Y|}, 
\eeq
 where $ \nu(\A)$ is given by \eqref{1.5}.
\end{definition}

Finally, we recall the next classical result of Campanato taken from \cite{C0} which is needed for the proof of our main result Theorem \ref{th2}:

\bt[Campanato] \label{th3}  Let $\mathcal{F},\mathcal{A} : \mathfrak{X} \larrow X$ be two mappings from the set $\mathfrak{X} \neq \emptyset$ into the Banach space $(X,\|\cdot\|)$. If there is a constant $K\in (0,1)$ such that
\beq \label{4.11}
\Big\|\mathcal{F}[u]-\mathcal{F}[v]-\big( \mathcal{A}[u]-\mathcal{A}[v]\big) \Big\| \, \leq\, K \big\| \mathcal{A}[u]-\mathcal{A}[v] \big\|
\eeq
for all $u,v \in \mathfrak{X}$ and if $\mathcal{A} : \mathfrak{X} \larrow X$ is a bijection, it follows that $\mathcal{F} :  \mathfrak{X} \larrow X$ is a bijection as well.
\et

\section{The AK-Condition of Ellipticity for Fully Nonlinear First Order Systems} \label{section3}
\ms

In this section we introduce and study a new ellipticity condition for the PDE system \eqref{1.1}which relaxes the K-Condition Definition \ref{def1} and still allows to prove existence and uniqueness of strong solutions to
\[
F(\cdot,\D u ) \,=\,  f, \ \ \ \text{ a.e.\ on }\R^n
\]
in the functional space \eqref{1.6}. Let 
\[
\A \ : \ \R^{Nn} \larrow \R^N
\]
be an elliptic reference linear map satisfying \eqref{2.1}.

\begin{definition}[The AK-Condition of ellipticity] \label{def2} Let  $n,N\geq2$  and 
\[
F\ :\ \ \R^{n} \by \R^{Nn} \larrow  \R^N
\]
a  Carath\'eodory map. We say  that $F$ is elliptic with respect to $\A$ when there exists a linear map
\[
\A \ : \ \R^{Nn} \larrow \R^N
\]
satisfying \eqref{1.3},  a positive function $\al $ with $\al, {1}/{\al} \in  L^{\infty}(\R^n)$   and $\be,\ga >0 $ with  $\be+\ga <1 $ such that
\beq \label{3.1}
\Big|\al(x) \Big[F(x,X+Y) -F(x,Y)\Big] \,-\,\A:X\Big| \, \leq \, \be \,\nu(\A) |X|\, +\, \ga\,| \A:X |.
\eeq
for all  $X,Y \in \R^{Nn} $ and a.e.\ $x\in \R^n$. Here $\nu(\A)$ is the ellipticity constant of  $A $ given by \eqref{1.5}.
\end{definition}

Nontrivial fully nonlinear examples of maps $F$ which are elliptic in the sense of the Definition \ref{def2} above are easy to find. Consider any fixed map $\A : \R^{Nn}\larrow \R^N$ for which $\nu(\A)>0$ and any Carath\'eodory map
\[
L\ : \ \R^n \by \R^{Nn}\larrow \R^N
\]
which is Lipschitz with respect to the second variable and
\[
\big\|L(x,\cdot)\big\|_{C^{0,1}(\R^{Nn})} \, \leq\, \be \, \nu(\A), \ \ \text{for a.e. }x\in \R^n
\]
for some $0<\be<1$. Let also $\al$ be a positive essentially bounded function with $1/\al$ essentially bounded as well. Then, the map $F :\R^n \by \R^{Nn}\larrow \R^N$ given by
\[
F(x, X)\, :=\, \frac{1}{\al(x)}\Big( \A:X\, +\, L(x, X) \Big)
\]
satisfies Definition \ref{def2}, since
\[
\begin{split}
\Big|\al(x)\Big[F(x, X+Y)-F(x,Y)\Big]-\A:X\Big|\, &\leq\, \big|L(x,X+Y)-L(x,Y) \big|\\
& \leq\, \be \, \nu(\A)|X| \\
& \leq\, \be \, \nu(\A)|X| \, +\, \frac{1-\be}{2}|\A:X|.
\end{split}
\]
As a consequence, $F$ satisfies the AK-Condition for the same function $\al(\cdot)$ and for the constants $\be$ and $\ga=(1-\be)/2$.

\ms

The following example shows that, given a reference tensor $\A$, there exist even \textbf{linear constant} \emph{``coefficients" $F$} which are elliptic with respect to $\A$ in the sense of our AK-Condition Definition \ref{def2} but which are \textbf{not} elliptic with respect to $\A$ in the sense of Definition \ref{def1} of \cite{K1}.

\begin{example} \label{ex9} Fix a constant $\al \in (0,1/2]$ and consider the linear map $F$ given by
\[
F(x, X)\, :=\,  \frac{1}{\al} \A:X, 
\]
where $A$ is the Cauchy-Riemann tensor of \eqref{CR}. Then, $F$ is elliptic in the sense of Definition \ref{def2} with respect to $\A$ for $\al(\cdot)\equiv \al$ and any $\be,\ga>0$ with $\be+\ga<1$, but it is not elliptic with respect to $\A$ in the sense of Definition \ref{def1}. Indeed for any $ X,Y \in \R^{Nn}$ we have:
\[
\begin{split}
\Big|\al \Big[F(\cdot,X+Y) -F(\cdot,Y)\Big] -\A:X \Big| \, 
&= \, \left|\al \left[\frac{1}{\al} \A:(X+Y) -\frac{1}{\al} \A:Y \right] -\A:X \right|   \\
&= \, 0 \\
&\leq \,\be \nu(\A) {|X|}\, +\, \ga\, | \A:X  |.
\end{split}
\]
On the other hand, by \eqref{CR} and \eqref{1.5} we have that $\nu(\A)=1$. Moreover, for 
\[
X_{0} \, :=\, 
\left[
\begin{array}{rr}
1 & 1       \\
1 & 1      \\
\end{array}
\right]
\]
we have $ |X_{0}  |=2$ and  $|\A:X_{0}| = 2$. Hence, for any $Y\in \R^{Nn}$ we have
\begin{align} \nonumber 
\Big|\Big[F(\cdot,X_{0}+Y) -F(\cdot,Y)\Big] -\A:X_{0} \Big| \, \nonumber 
&= \, \left| \left[\frac{1}{\al} \A:(X_{0}+Y) -\frac{1}{\al} \A:Y\right] -\A:X_{0} \right| \nonumber  \\
&=\, \left|\frac{1}{\al}\A:X_{0}\ -\A:X_{0} \right| \nonumber \\
&=\, \big|\A:X_{0} \big| \left|\frac{1}{\al}\ -1 \right|  \nonumber \\
&=\, 2 \left(\frac{1}{\al}\ -1 \right)  \nonumber \\
&\geq\, 2   \nonumber \\
&= \, \nu (A) \, {|X_{0}|} ,\nonumber
\end{align}
where we have used that $(1/\al)-1\geq 1$.  Our claim ensues.
\end{example}  

The essential point in the above example that makes Definition \ref{def2} more general than Definition \ref{def1} was the {introduction of the \emph{rescaling function $\al(\cdot)$.} Now we give a more elaborate example which shows that \emph{even if we ignore the rescaling function $\al$ and normalise it to $\al(\cdot)\equiv 1$, Definition \ref{def2} is still more general that Definition \ref{def1} with respect to the same fixed reference tensor $\A$.}

\begin{example} \label{ex8} Fix $c,b>0$ such that $c+b<1$ and $\sqrt{2} c+b>1$ and a unit vector $\eta\in \R^N$. Consider the Lipschitz function $F \in C^0 \big( \R^{2 \by 2} \big)$, given by: 
\beq \label{4.2}
F(x,X) := \A:X \, +\, \eta \cdot \big(\,b \, \big|X \big|+\,c \,\big| \A:X \big| \big),
\eeq
where $\A$ is again the Cauchy-Riemann tensor \eqref{CR}. Then, this $F$ satisfies
\beq \label{4.1}
\Big|\Big[F(\cdot,X+Y) -F(\cdot,X) \Big] -\A:Y\Big| \, \leq \, \be \,\nu(\A) {|Y|}\, +\, \ga\,\big| \A:Y \big|,
\eeq
for some  $\be,\ga >0 $ with  $\be+\ga <1$, but does not satisfy \eqref{4.1} with $\ga = 0$ \textbf{for any $0<\be<1$} for the same $\A$. Hence, $F$ satisfies Definition \ref{def2} (even if we fix $\al(\cdot)\equiv 1$) but it does not satisfy Definition \ref{def1}. Indeed we have:
\[
\begin{split} 
\Big|\A:Y - &\Big[F(\cdot,X+Y) -F(\cdot,X) \Big]\Big| \\
&= \, \Big|\A:Y\, -\, \A:Y \,-\, b  \eta  \Big(|X+Y| -|X| \Big) \,-\,  c  \eta  \Big( \big|\A:(X+Y)\big| - |\A:X | \Big) \Big| \\
&\leq\,  b | \eta| \Big| |X+Y| -|X| \Big| \,+ \,  c   | \eta| \Big| |\A:X+\A:Y| -|\A:X| \Big|   \\
&\leq\,  b | Y |\, +\,  c    |\A:Y |   
\end{split}
\]
and hence \eqref{4.1} holds for $\be = b$ and $\ga = c$. On the other hand, we choose
\[
X_0\,:=\,0,\ \ \ \ Y_0\, :=\, \, 
\left[
\begin{array}{cc}
1 & \zeta       \\
\zeta & 1      \\
\end{array}
\right],
\ \ \ \ \zeta\, := \, \frac{ 1-b  }{\sqrt{2 c^{2}- (1-b )^2}}. 
\]
This choice of $\zeta$ is admissible because our assumption $\sqrt{2} c+b>1$ implies $2 c^{2}- (1-b )^2 >0$. For these choices of  $X$ and $Y$, we calculate: 
\begin{align} 
\nonumber 
\left|\A:Y_0 - \Big[F(\cdot,X_0+Y_0) -F(\cdot,X_0) \Big]\right| \, \nonumber 
&=\,  \Big|\A:Y_0\, -\, F(\cdot,Y_0) \Big| \nonumber  \\
&=\,  \left|\A:Y_0\, -\, \A:Y_0\, - \, \eta  \Big(b \,|Y_0 | \,+\, c  |\A:Y_0 | \Big)\right| \nonumber \\
&=\,  | \eta| \Big|b|Y_0 | \,+ \, c  |\A:Y_0| \Big| \nonumber \\
&=\,  b \,| Y_0 | \ +  \ c  \, |\A:Y_0  | \nonumber .
\end{align}
We now show that 
\[
 b  \, | Y_0 | \, + \,  c |\A:Y_0 | \, =\, |Y_0|
\]
and this will allow us to conclude that \eqref{4.1} can not hold for any $\be<1$ if we impose $\ga=0$. Indeed, since $|Y_0|^2=2+2\zeta^2$ and $|\A:Y_0|^2=4\zeta^2$, we have
\[
\begin{split} 
\big(1-b \big)^2  | Y_0 |^2\, - \, c^2   |\A:Y_0  |^{2}  \, &=\, \big(1-b \big)^2 2\big(1+ \zeta ^{2}\big) \,-\, c^2  \, 4\zeta^{2}  \\
&=\,  2 \big(1-b \big)^{2} + 2 \Big(\big(1-b \big)^{2 }-2  c^2 \Big) \zeta^{2}  \\
&=\, 2 \big(1-b \big)^{2} + 2 \Big(\big(1-b \big)^{2 }-2  c^2 \Big)\frac{ (1-b)^2  }{2 c^{2}- (1-b )^2}\\
&=\, 0. 
\end{split}
\]
\end{example}

We now show that our ellipticity assumption implies a condition of pseudo-monotonicity coupled by a global Lipschitz continuity property. The statement and the proof are modelled after a similar result appearing in \cite{K2} which however was in the second order case.

\begin{lemma}[AK-Condition of ellipticity as Pseudo-Monotonicity] \label{pr1} Definition \ref{def2} implies the following statements:
\smallskip

There exist $\la>\ka>0$, a linear map $\A : \R^{Nn}\larrow \R^N$ satisfying \eqref{1.3} a positive function $\al$ such that  $\al,1/\al \in L^{\infty}(\R^n)$ with respect to which $ F$ satisfies 
\beq \label{4.3}
 (\A:Y)^\top\Big[F(x,X+Y) - F(x,X) \Big]\ \geq  \ \frac{ \la}{ \al(x)}|\A:Y|^2\, -\, \frac{\ka}{ \al(x)}\nu(\A)^2|Y|^2,
\eeq
for all  $X,Y \in \R^{Nn} $ and a.e.\ $x \in \R^n$. In addition, $F(x,\cdot)$ is  Lipschitz continuous on $\R^{Nn}$,  essentially uniformly in $x\in \R^n$; namely, there exists $M>0$ such that 
\beq \label{4.4}
 \big| F(x,X) \,-\, F(x,Y)\big| \, \leq\, M |X-Y|
\eeq
for a.e.\ $x\in \R^n$ and all $X,Y\in \R^{Nn}$.
\end{lemma}

\BPL \ref{pr1}.  Suppose that  Definition \ref{def2} holds for some constant $\be,\ga>0$ with $\be+\ga<1$, some positive function $\al$ with $\al,1/\al \in L^{\infty}(\R^n)$ and some linear map $\A : \R^{Nn}\larrow \R^N$ satisfying \eqref{1.3}. Fix $\e>0$. Then, for a.e.\ $x\in \R^N$ and all $X,Y\in \R^{Nn}$ we have:
\[
\begin{split}  
&| \A: Y |^2\, + \,  \al(x)^2 \Big|F(x,X+Y) -F(x,X) \Big|^2 \\
&- \,2\al(x)  \left( \A: Y \right)^\top\Big[  F(x,X+Y) -F(x,X) \Big] \\   
&\ \  \ \ \ \ \ \ \leq \phantom{\Big|}
\, \be^2 \nu(\A)^2|Y|^2 \, + \, \ga^2 | \A: Y |^2\,  + \, 2  \be  \nu(\A) |Y| \, \ga | \A: Y |  
\end{split}
\]
which implies
\[
\begin{split}  
| \A: Y |^2\,-& \, 2\al(x)  \left( \A: Y \right)^\top\Big[  F(x,X+Y) -F(x,X) \Big] \\   
&\leq \,\be^2 \nu(\A)^2|Y|^2 + \ga^2 | \A: Y |^2 \,+\, \frac{\be^2  \nu(\A)^2 |Y|^2  }{ \e} \,+\, \e \ga^2  | \A: Y |^2  .
\end{split}
\]
Hence,
\[
\begin{split}  
 \left( \A: Y \right)^\top & \Big[  F(x,X+Y) -F(x,X) \Big] \\
&\geq \, \frac{1}{\al(x) }\left( \frac{1- \ga^2 -  \e \ga^2}{2} \right) | \A: Y |^2 \, -\,   \frac{1}{\al(x) }\left(\frac{\e \be^2 +\be^2 }{2 \e} \right)  \nu(\A)^2|Y|^2.
\end{split}
\]
By choosing $\e :={\be}/{\ga}$, from the above inequality we obtain \eqref{4.3} for the values
\[
\la\, :=\, \frac{1-\ga(\ga+\be)}{2},\ \ \ \ \ka\,:=\,\frac{\be(\ga+\be)}{2}.
\]
These are admissible because $\ka>0$ and $\la >\ka$ since 
\[
\la \,-\,\ka \,=\, \frac{1-(\be+\ga)^2 }{2}\, >\, 0.
\]
In addition, again by  \eqref{3.1} we have:
\beq
\al(x) \Big|F(x,X) -F(x,Y)\Big|  \, \leq  \,  \be \nu(\A) { |X-Y|}\,+\,\ga\,\big|\A:(X-Y)\big | \,+\,\big|\A:(X-Y)\big | ,  \nonumber
\eeq
and hence,
\[
\begin{split}  
\Big|F(x,X) -F(x,Y)\Big|  \,
&\leq  \,\frac{1}{\al(x) } \Big( (1+\ga) \big|\A:(X-Y)\big | \,+\, \be \nu(\A) \big|X-Y \big|\Big) \\
&\leq \, \Bigg\{\left\|\frac{1}{\al(\cdot) } \right\|_{ L^{\infty}(\R^n)}  \Big( (1+\ga)  |\A | + \be \nu(\A)\Big)\Bigg\} |X -Y| \nonumber
\end{split}
\]
 for a.e.\ $x\in \R^N$ and all $X,Y\in \R^{Nn}$, which immediately leads to \eqref{4.4} and the proposition ensues.       \qed
\ms

\section{Well-Posedness of Global Fully Nonlinear First Order Elliptic Systems} \label{section4}

In this section we state and prove the main result of this paper which is the following:

\bt[Existence-Uniqueness] \label{th2} Assume that $n\geq 3$, $N\geq 2$ and let $F : \R^n \by \R^{N \by n}\larrow \R^N$ be a Carath\'eodory map, satisfying Definition \ref{def2} with respect to a reference tensor $\A$ which satisfies \eqref{2.1}. 

\ms

\noi (1) For any two maps $v,u \in W^{1;2^*\!,2}(\R^n,\R^N)$ (see \eqref{1.6}), we have the estimate 
\beq \label{4.10}
\|v-u\|_{W^{1;2^*\!,2}(\R^n)}\, \leq\, C \big\|F(\cdot,Dv)-F(\cdot,Du) \big\|_{L^{2}(\R^n)}
\eeq
for some $C>0$ depending only on $F$. Hence, the PDE system $F(\cdot,\D u) = f$ has at most one solution.

\ms 

\noi (2) Suppose further that $F(x,0)=0$ for a.e.\ $x\in \R^n$. Then for any $f\in L^2(\R^n,\R^N)$, the system
\[
F(\cdot,\D u)\, =\, f, \ \ \text{ a.e.\ on }\R^n,
\]
has a unique solution $u$ in the space $W^{1;2^*\!,2}(\R^n,\R^N)$ which also satisfies the estimate
\beq \label{4.9}
\|u\|_{W^{1;2^*\!,2}(\R^n)}\, \leq\, C \|f\|_{L^{2}(\R^n)}
\eeq
for some $C>0$ depending only on $F$. 
\et

\BPT \ref{th2}. (1) Let $\al$ and $\A$ be as in Definition \ref{def2} and fix $u,v\in W^{1;2^*\!,2}(\R^n,\R^N)$. Since $\A$ satisfies \eqref{2.1}, by Plancherel's theorem (see e.g.\ \cite{F}) we have:
\begin{align} \label{4.14}
\frac{1}{\nu(\A)}\big\|\A:\big(\D v -\D u\big) \big\|_{L^2(\R^n)}\, &= \, \frac{1}{\nu(\A)}\big\|\A:\big(\widehat{Dv} - \widehat{\D u}\big) \big\|_{L^2(\R^n)} \nonumber\\
&= \, \frac{1}{\nu(\A)}\big\| \A: \big(\widehat{v}-\widehat{u}\big) \ot (2\pi i \mathrm{Id}) \big\|_{L^2(\R^n)}   \nonumber \\
&\geq \,  \big\| \big(\widehat{v}-\widehat{u}\big) \ot (2\pi i \mathrm{Id}) \big\|_{L^2(\R^n)} \\
&= \,  \big\|\widehat{\D v} - \widehat{\D u} \big\|_{L^2(\R^n)}  \nonumber\\
&= \,  \big\| {\D v} -  {\D u} \big\|_{L^2(\R^n)},  \nonumber
\end{align}
where we symbolised the identity map by ``Id", which means $\mathrm{Id}(x):=x$. Further, by Definition \ref{def2} also we have
\begin{align} \label{4.15}
\Big\|\al(\cdot) \big[ & F(\cdot,\D u) -  F(\cdot,\D v)\big] \,-\, \A: \big(\D u-\D v \big)\Big\|_{L^2(\R^n)}\nonumber \\
&\leq \, \be \nu(\A)\big\| {\D u} -  \D v \big\|_{L^2(\R^n)}+\, \ga  \big\|\A: \big(\D u-\D v \big)  \big\|_{L^2(\R^n)} \nonumber
\end{align}
Using the estimate \eqref{4.14} above this gives:
\begin{align} 
\Big\|\al(\cdot) \big[ & F(\cdot,\D u) -  F(\cdot,Dv)\big] \,-\, \A: \big(\D u-\D v \big)\Big\|_{L^2(\R^n)} \nonumber\\
&\leq \, \be \big\|\A: \big(\D u-\D v \big)  \big\|_{L^2(\R^n)}+ \, \ga  \big\|\A: \big(\D u-\D v \big)  \big\|_{L^2(\R^n)}\\
&\leq \, \big(\be +\ga\big)  \big\|\A: \big(\D u-\D v \big)  \big\|_{L^2(\R^n)}
\nonumber
\end{align}
and hence
\begin{align} 
\nonumber
\big(\be +\ga\big)  & \big\|\A: \big(\D u-\D v \big)  \big\|_{L^2(\R^n)} \nonumber\\
&\geq \, \Big\|\A: \big(\D u-Dv \big) \,-\, \al(\cdot) \big[F(\cdot,\D u) -  F(\cdot,\D v)\big]\Big\|_{L^2(\R^n)}\nonumber\\
&\geq \, \big\|\A: \big(\D u-\D v \big)\big\|_{L^2(\R^n)} -\, \Big\|\al(\cdot) \big[F(\cdot,\D u) -  F(\cdot,\D v)\big]\Big\|_{L^2(\R^n)}\nonumber
\end{align}
which implies the following estimate:
\begin{align} \nonumber
\Big\|\al(\cdot) \big[F(\cdot,\D u) -  F(\cdot,\D v)\big]\Big\|_{L^2(\R^n)}\nonumber
&\geq\big[1-(\be +\ga)\big]\big\|\A: \big(\D u-\D v \big)]\big\|_{L^2(\R^n)}\nonumber\\
&\geq\big[1-(\be +\ga)\big] \nu(\A)\big\|\D u-\D v\big\|_{L^2(\R^n)}\nonumber
\end{align}
Since $\be +\ga<1$, we have the estimate:
\begin{align} \label{4.16}
\frac{\|\al(\cdot) \|_{L^\infty(\R^n)}}{\big[1-(\be +\ga)\big] \nu(\A)} \big\|F(\cdot,\D u) -  F(\cdot,\D v)\big\|_{L^2(\R^n)}\, 
&\geq\, \big\| \D u-\D v \big\|_{L^2(\R^n)}.
\end{align}
By  \eqref{4.16}, and the fact that $n \geq 3$, the Gagliardo-Nirenberg-Sobolev inequality gives the estimate
\begin{align}
\| u -v\|_{W^{1;2^*\!,2}(\R^n)} \, \leq\, C \big\| F(\cdot,\D u) -  F(\cdot,\D v) \big\|_{L^2(\R^n)}
\end{align}
where $C>0$ depends only on $F$.

\ms

\noi (2) By our assumptions on $F$ and that $F(x,0)=0$, Lemma \ref{pr1} implies that there exists an $M>0$ depending only on $F$, such that for any $u\in W^{1;2^*\!,2}(\R^n,\R^N)$, we have the estimates
\begin{align} \label{4.12}
\big\|\al(\cdot)F(\cdot,\D u) \big\|_{L^2(\R^n)}
&\,=\Big\|\al(\cdot)\big[F(\cdot,0+\D u)-F(\cdot,0)\big] \Big\|_{L^2(\R^n)}  \nonumber\\
&=\, M  \| \al(\cdot)\|_{L^{\infty}(\R^n)}\| \D u\|_{L^2(\R^n)} \\
&\leq\, M \| \al(\cdot)\|_{L^{\infty}(\R^n)}\| u\|_{W^{1;2^*\!,2}(\R^n)} \nonumber
\end{align}
and also
\beq \label{4.13}
\|\A:\D u\|_{L^2(\R^n)}\, \leq\, \|\A\|\, \| \D u\|_{L^2(\R^n)} \, \leq\,  \|\A\| \| u\|_{W^{1;2^*\!,2}(\R^n)}.
\eeq
We conclude from \eqref{4.12} and \eqref{4.13}  that the differential operators 
\[
\left\{
\begin{array}{l}
\mathscr{A}[u]\ :=\ \A :\D u, \ms \ms\\
\mathscr{F}[u]\ :=\ \al(\cdot)F(\cdot, \D u),
\end{array}
\right.
\]
map the functional space $W^{1;2^*\!,2}(\R^n,\R^N)$ into the space $L^2(\R^n,\R^N)$. Note that Theorem \ref{th1} proved in \cite{K1} implies that the linear operator
\[
A\ :\  W^{1;2^*\!,2}(\R^n,\R^N) \larrow  L^2(\R^n,\R^N) 
\]
is a bijection. Hence, in view of inequality \eqref{4.15} above and the fact that $\be +\ga <1$, Campanato's nearness Theorem \ref{th3} implies that $F$ is a bijection as well. As a result, for any $g\in L^2(\R^n,\R^N)$, the PDE system 
\[
\al(\cdot) F(\cdot, \D u)\, =\, g,\ \  \text{ a.e. on }\R^n,
\]
has a unique solution $u\in W^{1;2^*\!,2}(\R^n,\R^N)$. Since $\al(\cdot),{1}/{\al(\cdot)} \in  L^{\infty}(\R^n)$, by selecting $g=\al(\cdot) f$, we conclude that the problem
\[
 F(\cdot, \D u)\, =\, f,\ \  \text{ a.e. on }\R^n,
\]
has a unique solution in $ W^{1;2^*\!,2}(\R^n,\R^N)$. The proof of the theorem is now complete.      \qed
\ms 

\ms

\bibliographystyle{amsplain}

\end{document}